\documentclass[a4paper,11pt]{article}
\usepackage{amssymb}
\usepackage{macros, margins, mysect}
\usepackage{graphicx}
\usepackage[font=small]{caption}
\usepackage{framed}

\begin{document}
\title{Runoff on rooted trees}
\author{Owen Dafydd Jones}
\date{July 2018}
\maketitle

\abstract{
We introduce an idealised model for overland flow generated by rain falling on a hill-slope.
Our prime motivation is to show how the coalescence of runoff streams promotes the total generation of runoff.
We show that, for our model, as the rate of rainfall increases in relation to the soil infiltration rate, there is a distinct phase-change.
For low rainfall (the subcritical case) only the bottom of the hill-slope contributes to the total overland runoff, while for high rainfall (the supercritical case) the whole slope contributes and the total runoff increases dramatically.
We identify the critical point at which the phase-change occurs, and show how it depends on the degree of coalescence.
When there is no stream coalescence the critical point occurs when the rainfall rate equals the average infiltration rate, but when we allow coalescence the critical point occurs when the rainfall rate is less than the average infiltration rate, and increasing the amount of coalescence increases the total expected runoff.
}

\section{Context and main results}

The motivation for this work is the problem of modelling surface runoff.
Surface runoff depends on rainfall, infiltration into the soil, and surface topography, all of which vary spatially and temporally.
In particular, spatial variation of infiltration and topographical variability make it difficult to fit differential fluid flow models (based on the Navier-Stokes equations) at coarse scales, because parameters lose their physical meaning \cite{GMM92}, while fitting them at fine scales requires high-resolution data and is numerically prohibitively slow.

A practical alternative to measuring the infiltration and topography of a hill-slope at high-resolution is to summarise small-scale variation statistically, which leads naturally to stochastic runoff models.
Suppose that we divide our hill-slope into cells.
We can model spatial variation of the soil infiltration by supposing that for each cell it is sampled independently from some infiltration distribution.
It is known that this variation is enough for a hill-slope to produce surface runoff even when the rainfall is less than the average infiltration \cite{JSL13, JLS16}.
In this work we consider in addition the effect of variation in the micro-topography from one cell to the next.

When surface runoff forms on a hill-slope, we see small trickles combining to form larger rivulets, which are proportionally less susceptible to soil infiltration.
We call this mechanism {\em coalescence}, and we want to show how it impacts surface runoff.
We suppose that micro-topography will affect the direction of runoff from a cell.
Water necessarily flows downhill, but local variation in the topography can mean that instead of taking the most direct route down a slope, a rivulet is diverted to the left or right.
We can model this by adding randomness to the direction in which runoff flows out of a cell, where the degree of randomness reflects the roughness of the hill-slope.

\subsection{An illustrative simulation}

We can explore the effect of coalescence with a simple simulation.

Divide a hill-slope into a rectangular $m \times n$ grid of cells, so that cells $(1,j)$ are at the top of the slope and cells $(m,j)$ at the bottom.
We suppose that if there is any runoff from cell $(i,j)$, then it can run to cell $(i+1,j-1)$, $(i+1,j)$ or $(i+1,j+1)$ with probabilities $\de$, $1-2\de$ and $\de$ respectively.
This direction does not change over time, and we don't allow runoff to exit via the sides of the grid.

Next suppose that rainfall is constant rate $\rho$, and that the infiltration rate in cell $(i,j)$ has an exponential distribution with mean $1$, independent of other cells.
We suppose that the system is in temporal equilibrium, so that at each cell the rates of rainfall, infiltration, runon from above, and runoff do not vary in time.
Let $J_{(i,j)}$ be the infiltration rate for cell $(i,j)$ and $W_{(i,j)}$ be the runoff rate from cell $(i,j)$, then we have
\begin{equation}\label{eqn1}
W_{(i,j)} = \left( \rho - J_{(i,j)} + \sum_{k:(i-1,k) \to (i,j)} W_{(i-1,k)} \right) \vee 0,
\end{equation}
where we write $(i-1,k) \to (i,j)$ if runoff from $(i-1,k)$ runs on to cell $(i,j)$ (in which case $k \in \{j-1,j,j+1\}$).

\begin{figure}
\begin{center}
\includegraphics[width=14cm]{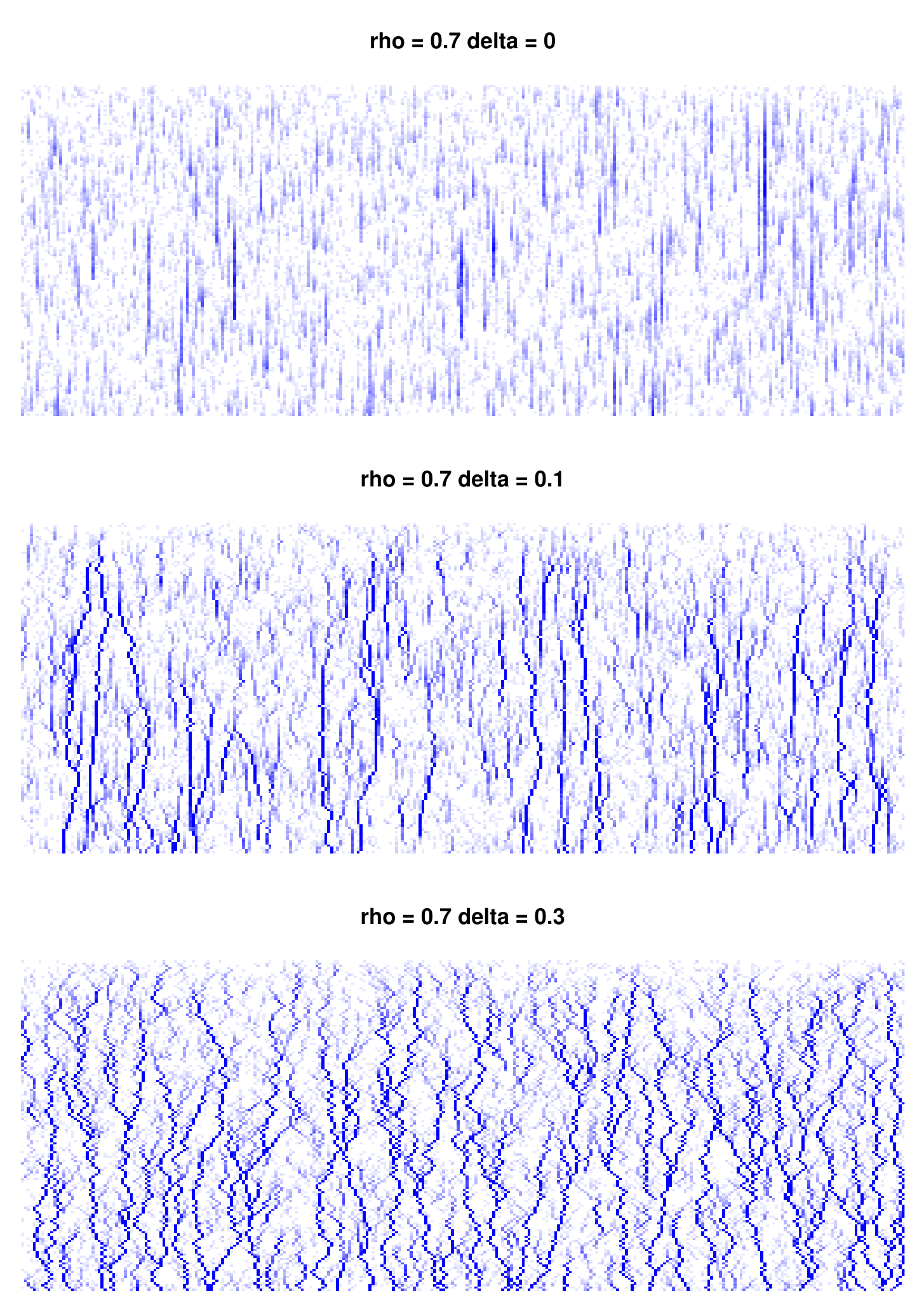}
\caption{Simulation output illustrating the effect of coalescing streams.
In each case water flows from top to bottom; the darker the pixel the greater the flow.
Each cell has rainfall rate $\rho$, a random infiltration rate with mean $1$, a chance $\de$ that runoff is directed down and to the left, and a chance $\de$ that it is directed down and to the right.
For $\de=0$ any runoff that is formed is eventually reabsorbed back into the slope, but for $\de=0.3$ runoff makes its way down the whole slope, even though the average infiltration rate exceeds the rainfall rate. 
Note that the three plots have been scaled so that the maximum runoff is the same shade; the maximum runoff actually increases with $\de$.
The code can be found at {\tt http://researchers.ms.unimelb.edu.au/\~{}apro@unimelb/spuRs/index.html}.}
\label{sim.fig}
\end{center}
\end{figure}

In Figure \ref{sim.fig} we give three realisations of this process, for $m=150$, $n=300$, $\rho = 0.7$ and $\de \in \{0, 0.1, 0.3\}$.
We see that for $\de=0$ any runoff that is formed is eventually reabsorbed back into the slope, but for $\de=0.3$ runoff makes its way down the whole slope, even though the average infiltration rate exceeds the rainfall rate.
This pattern is consistent for different values of $\rho$---as $\de$ increases there is increasing runoff---moreover for any $\rho$ there is a threshold value of $\de$ after which we start to see runoff making its way down the whole hill-slope from top to bottom.
It what follows we will develop an abstract model for runoff for which we can establish this phase-change precisely.

\subsection{Drainage trees}

\begin{figure}
\begin{center}
\includegraphics[width=6cm]{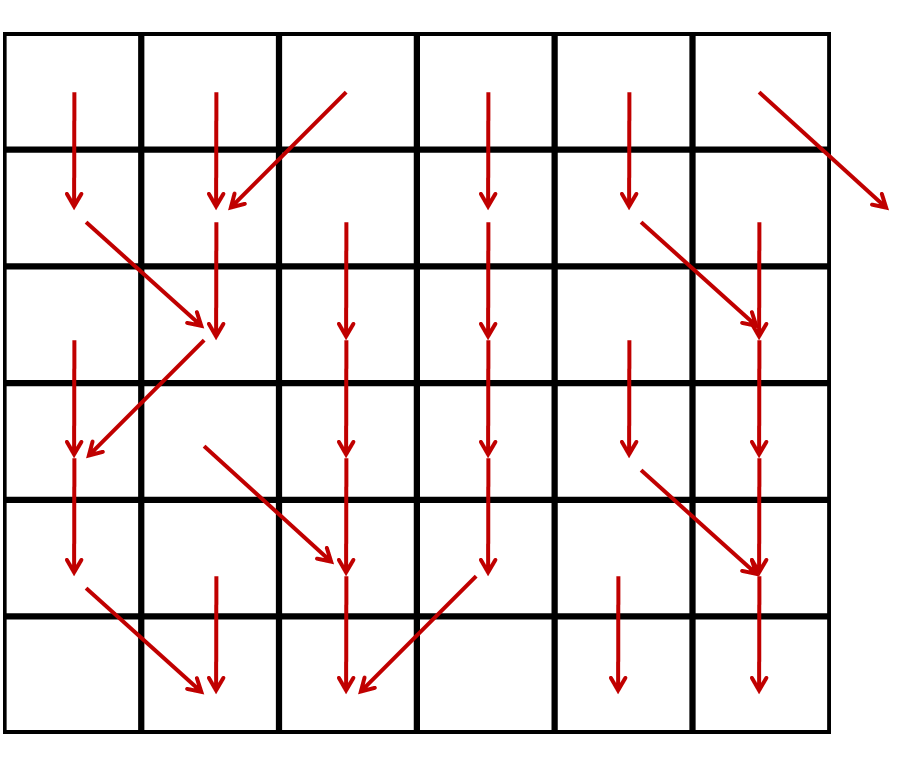}
\includegraphics[width=6cm]{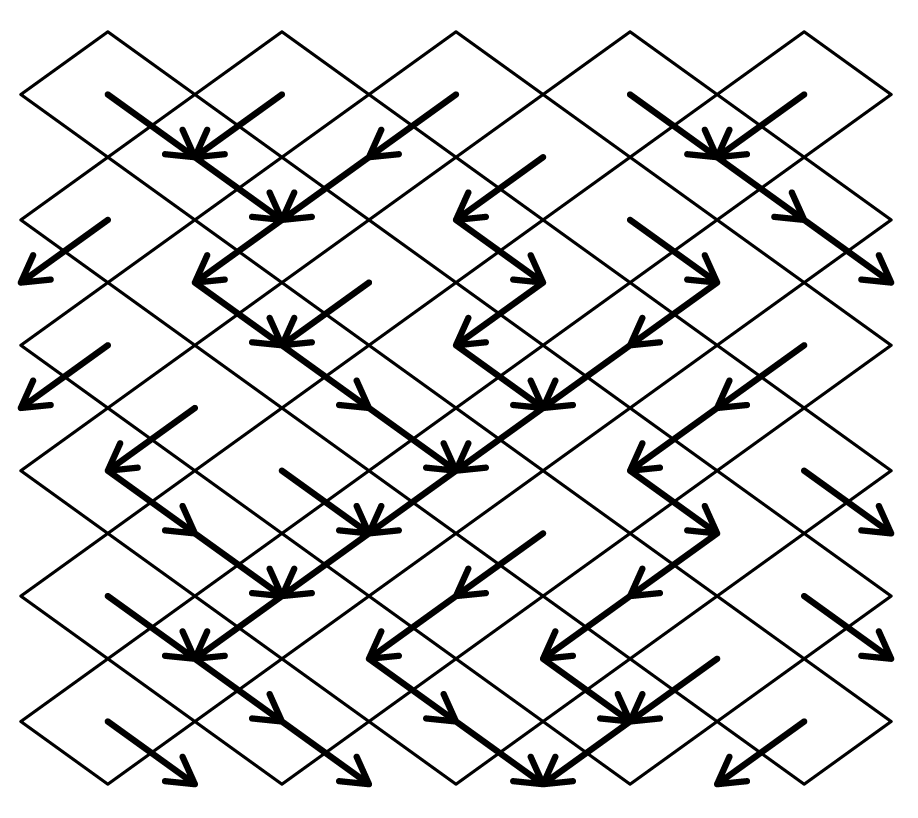}
\caption{Two realisations of paths of potential runoff, for a square lattice (left) and a diamond lattice (right).
In either case, if you take a single cell at the bottom and consider all the cells that could drain into it, you get a tree.}
\label{lattice.fig}
\end{center}
\end{figure}

Suppose that we have a hill-slope divided into cells, and that the runoff from any given cell will flow into a unique cell below.
For example, using the square lattice we could allow runoff into any one of the three cells below with a common edge or vertex, or using the diamond lattice might choose to allow runoff into either of the two cells below with a common edge.
Examples of the runoff paths you can generate in these two cases are given in Figure \ref{lattice.fig}.
Selecting a single cell at the bottom of the hill-slope and then considering all the cells that could potentially drain into it, we get a rooted tree, which we call a {\em drainage tree}, where the nodes correspond to cells and edges indicate where runoff flows from one cell to the next (runoff is always towards the root).

In what follows we will consider runoff on a single drainage tree.
Given a rooted tree, let $X_i$ be the difference between rainfall and infiltration in cell $i$, and denote by $\{j:j\to i\}$ those cells that send runoff directly into node $i$.
The runoff from cell/node $i$ is then
\begin{equation}\label{runoff.eqn}
W_i =\left(  X_i + \sum_{j:j\to i} W_j \right) \vee 0.
\end{equation}

We can view the $X_i$ as either rates or, if we integrate over discrete time periods, volumes.
In either case we are assuming that the system is temporally homogeneous and in equilibrium.

Our object of interest is $W_0$, where $0$ is the root.
If the tree is finite then $W_0$ is clearly well defined: we have $W_j = X_j \vee 0$ for all leaves $j$, then we can use Equation \ref{runoff.eqn} to recursively calculate $W_i$ for all other nodes.
For infinite trees we just define $W_0 = \lim_{n\to\infty} W_0^{(n)}$, where $W_0^{(n)}$ is the root-runoff from the tree truncated at generation/height $n$.
The $W_0^{(n)}$ are clearly non-decreasing so $W_0$ exists, though may be improper.

It is interesting to note that if instead of working our way down from the leaves we consider working our way up from the root, then we get
\begin{equation}\label{maxsum.eqn}
W_0 = 
\max_{T \in \{ \mbox{rooted subtrees} \}}
\sum_{i \in T} X_i
\end{equation}
where a rooted subtree is any subtree including the root $0$, or the empty subtree (in which case we take the sum to be $0$).

\subsection{Drainage trees from the diamond lattice}

In hydrology the use of the diamond lattice to generate random drainage patterns can be traced back to Scheidegger \cite{Sch67}.
Note however that in the hydrological literature drainage trees are used to describe river networks, rather than the small-scale patterns of ephemeral surface runoff that we are interested in.

Take as our hill-slope a half plane extending upwards ad infinitum, and divide it into cells using a diamond lattice.
Furthermore suppose that from each cell runoff goes left with probability $\be \in (0, 1)$ and right with probability $\bebar = 1 - \be$, independently of all other cells.
Consider the drainage tree attached to a single cell at the bottom of the slope.
If we think of the tree as growing from its root, then for any node the number of offspring has distribution

\begin{equation}\label{offspring.eqn}
\begin{array}{c|ccc}
z & 0 & 1 & 2 \\
\hline
\Pb(Z=z) & \be\bebar & \be^2 + \bebar^2 & \be\bebar.
\end{array}
\end{equation} 

If $\be = 0$ or $1$ then our tree degenerates to one-dimension, but is infinite in size.
This case is considered in \cite{JSL13, JLS16}, and from here on we will assume that $\be \in (0, 1/2]$, unless stated otherwise.
We can interpret $\be$ in terms of surface roughness: $\be=1/2$ gives the roughest surface, with decreasing values corresponding to smoother surfaces.
In terms of our model, we get the greatest degree of coalescence when $\be=1/2$, and the least (none) when $\be = 0$.
The cases $\be$ and $1-\be$ are equivalent by symmetry.

Clearly the offspring in any given generation are dependent: letting $Z_n$ be the number of nodes in generation $n$ (where the root is generation $0$), we have that $Z_{n+1} = Z_n + D_n - 1$, where $D_n$ is independent of $Z_n$ and is distributed as (\ref{offspring.eqn}).
That is, the tree diameter is given by a random walk with zero drift, from which it follows immediately that the tree is almost surely finite, but that its expected height is infinite.

In what follows we will approximate the diamond lattice tree with a critical Bienaym\'e-Galton-Watson (BGW) tree, by the simple expediency of dropping the dependence between offspring numbers in each generation.
That is, we use the offspring distribution (\ref{offspring.eqn}).

\subsection{Runoff on a critical BGW tree}

Suppose that we are given a critical BGW tree with offspring distribution (\ref{offspring.eqn}), for $\be \in (0, 1/2]$.
We associate i.i.d.\ random variables $X_i$ with each node, and are interested in the $W_i$ as defined by (\ref{runoff.eqn}).

Write $X$ and $W$ for $X_0$ and $W_0$, the point contribution and nett runoff at the root, and let $W_L$ be the runoff from the cell above left and $W_R$ the runoff from the cell above right.
From the self-similar structure of the BGW tree we have that $W$, $W_L$ and $W_R$ are identically distributed, and $W_L$ and $W_R$ are independent.

Let $I_L$ indicate if the cell above left drains into the root cell, and similarly for $I_R$, then $\Ex I_L = \bebar$, $\Ex I_R = \be$, and
\begin{equation} \label{W.eqn}
W = (W_L I_L + W_R I_R + X) \vee 0.
\end{equation}
This equation is the focus of the remainder of the paper, and we note to begin with that that an almost surely finite solution $W$ exists, because the tree is almost surely finite and we can construct solutions using (\ref{runoff.eqn}).

\begin{figure}
\begin{center}
\includegraphics[width=8cm]{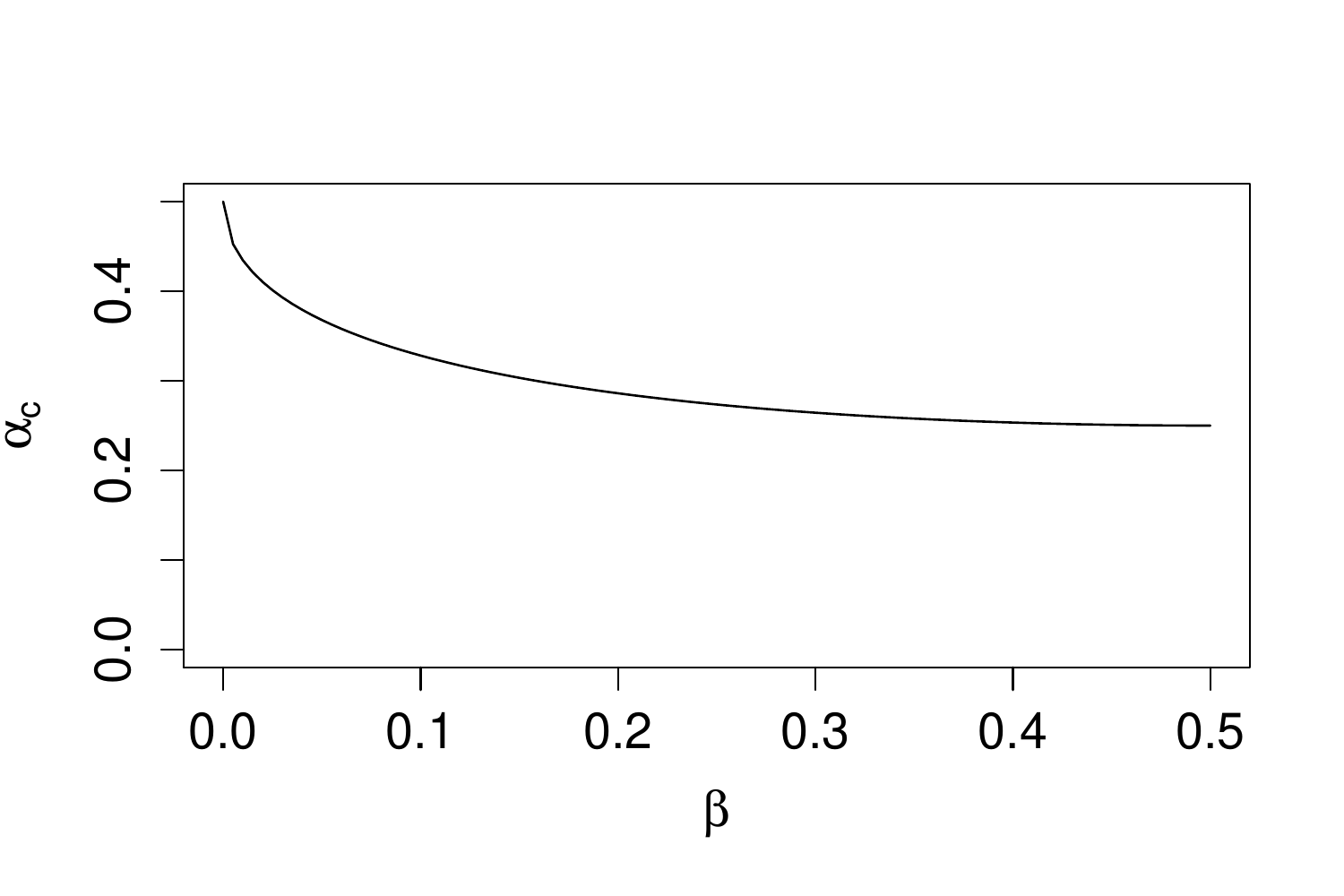}
\caption{The function $\al_c(\be) = (1/2) \left( 1 + \be\bebar - \sqrt{\be\bebar (2 + \be\bebar)} \right)$.}
\label{alpha-c.fig}
\end{center}
\end{figure}

We summarise our main results here, and defer the proofs to the next two sections.
Section \ref{left_cts.sec} then discusses some generalisations.
In order to obtain exact results, we will suppose that $X$ has the following distribution, for some $\al \in (0, 1)$,

\begin{equation}\label{X.eqn}
\begin{array}{c|cc}
x & 1 & -1 \\
\hline
\Pb(X=x) & \al & \albar = 1-\al
\end{array}
\end{equation}

If $\al = 0$ then $W = 0$, while if $\al = 1$ then $W$ is just the size of the tree.

\begin{prop}
Put
\begin{equation}\label{a_c.eqn}
\al_c(\be) = (1/2) \left( 1 + \be\bebar - \sqrt{\be\bebar (2 + \be\bebar)} \right)
\end{equation}
then for $\al \leq \al_c(\be)$
\begin{equation}\label{Ex_W.eqn}
\Ex W = \frac 1{2\be(1-\be)} \left( 1 - 2\al - \sqrt{1 - 4\al(1 - \al + \be(1-\be))} \right) 
\end{equation}
while for $\al > \al_c(\be)$, $\Ex W = \infty$.
\end{prop}

\begin{figure}
\begin{center}
\includegraphics[width=8cm]{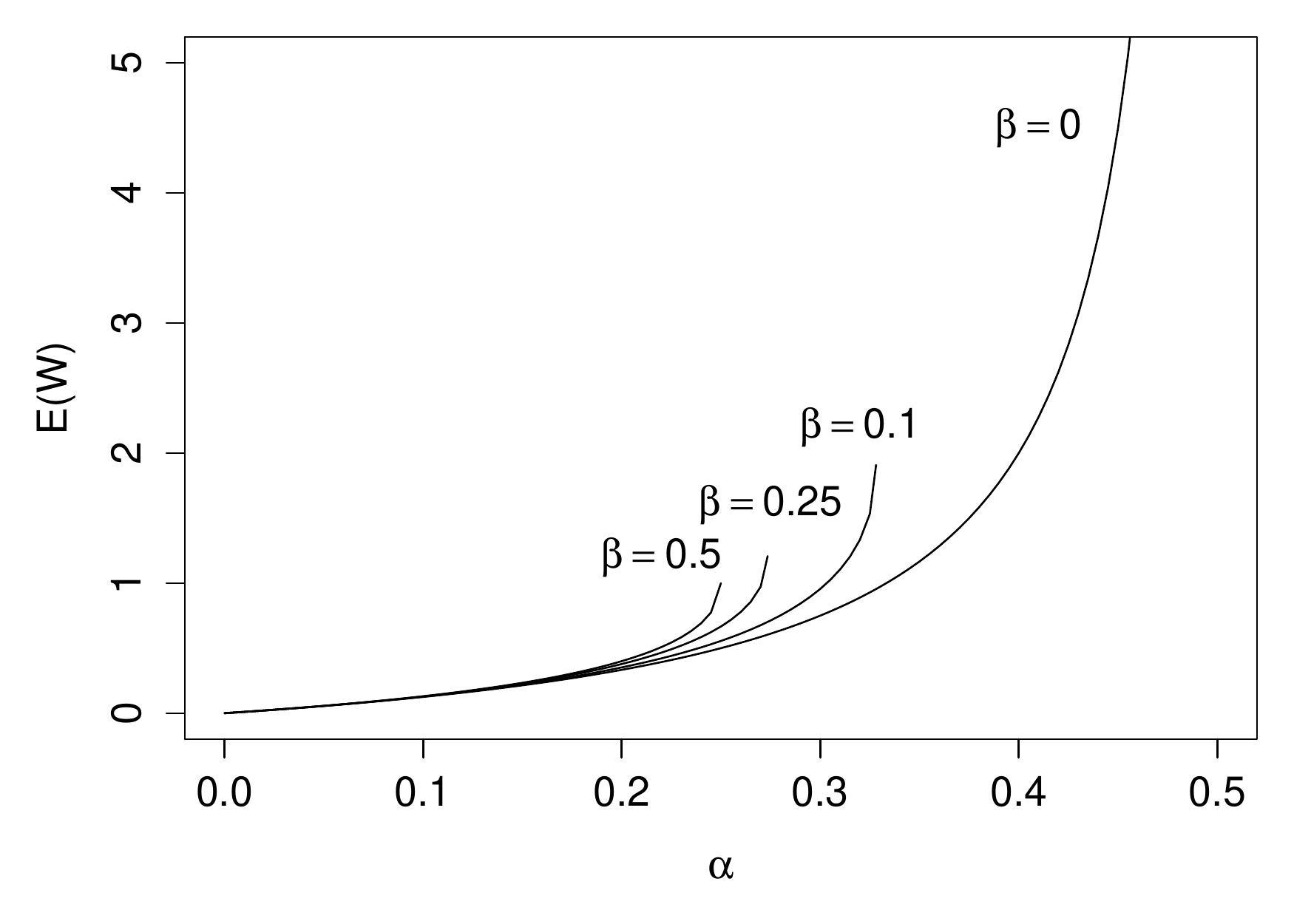}
\caption{Expected runoff for various values of $\al$ and $\be$.
For $\al$ larger than the given ranges, the expected mean is $\infty$.}
\label{ExW.fig}
\end{center}
\end{figure}

Figure \ref{alpha-c.fig} gives a plot of $\al_c$.
Figure \ref{ExW.fig} gives plots of $\Ex W$ against $\al$ for various $\be$.

Note that for all $\be > 0$ we have $\al_c(\be) < 1/2$, so that $\Ex X < 0$.
That is, the critical point at which the expected runoff becomes infinite happens when the rainfall is less than the expected infiltration.
Moreover, as the degree of coalescence increases, that is as $\be \upto 1/2$, less rainfall is required to reach the critical point.
We say that the runoff is subcritical/critical/supercritical as $\al$ is less than/equal to/greater than $\al_c$.

We can say more about the size of $W$, in terms of how heavy its tail is.

\begin{prop}
Put
\begin{equation}\label{h.eqn}
h(t) = \frac {t[1 - \al(4\be\bebar + t) - \al^2(1-2\be)^2(1-t^2)]}
{4\albar\be^2\bebar^2(1 - \al(1-t^2))}
\end{equation}
and let $t_0$ be the point at which $h$ achieves its maximum in $[0, 1]$,
then, as $x \to \infty$,
\begin{equation}\label{main.eqn}
\Pb(W > x) \sim
\left\{\begin{array}{ll}
\sqrt{ \frac {-h''(1)(1-\al)}{8\pi}} x^{-3/2} & \al = \al_c \\
\sqrt{ \frac {(h(t_0) - h(1))(1-\al)}\pi } x^{-1/2} & \al > \al_c
\end{array}\right.
\end{equation}
For $\al < \al_c$, $W$ has all positive moments finite.
\end{prop}

Let $N_T$ be the total number of nodes in our drainage tree, then it is known that
\begin{equation}\label{NT.eqn}
\Pb(N_T = n) \approx n^{-3/2} \frac 1{2\sqrt{\pi\be\bebar}}.
\end{equation}
We give a sketch of the proof in Section \ref{size_tree.sec}.
In particular, even though $N_T$ is almost surely finite, we have $\Ex N_T = \infty$.
This suggests that when $\Ex W = \infty$ it is because the runoff at the root is getting contributions from all over the tree, where we say that a node $i$ {\em contributes} to the runoff at the root if there is a path $(i_1 = i, i_2, \ldots, i_n = 0)$ from $i$ to the root such that $W_{i_k} > 0$ for all $k$.
Our final result for this section says that this is indeed the case.
Moreover, we see that when $\Ex W < \infty$ only the bottom of the tree is contributing.

We can re-express (\ref{W.eqn}) as
\begin{eqnarray}
W &=& (W_LI_L + W_RI_R + X) \vee 0 \nonumber \\
&=& W_LI_L + W_RI_R + Y \label{Y.eqn}
\end{eqnarray}
where $Y$ is the nett contribution to runoff from our cell.
$W_L$, $W_R$, $I_L$, $I_R$ and $X$ are all independent, but the distribution of $Y$ depends on $X$, $W_L$, $W_R$, $I_L$ and $I_R$.

\begin{prop}\label{contribution.prop}
When $\al > \al_c$ we have $\Ex Y > 0$, and for every $\de > 0$ there is an $\ep > 0$ such that
\[
\Pb( \mbox{$100(1-\de)\%$ of the tree height contributes to root runoff} ) \geq \ep.
\]
When $\al \leq \al_c$ we have $\Ex Y = 0$, and the expected tree height contributing to root runoff is finite. 
\end{prop}

We finish the section with a hydrologically inclined qualitative summary of our results.

\begin{framed}
When the rainfall is subcritical, only the bottom of the hill-slope contributes to runoff, but when the rainfall is supercritical, the whole slope starts contributing, giving a phase-change in the amount of runoff produced.
The critical point depends on the degree of coalescence induced by the micro-topography of the hillslope.
\end{framed}

\subsection{Links to other work}

The pattern of drainage trees produced by the diamond lattice is a familiar probabilistic object, known as {\em coalescing random walks} or the {\em voter model} in dimension one.
See for example the books of Liggett \cite{Lig85} or Durrett \cite{Dur88}.

As noted above, the offspring numbers in generation $n$ of a drainage tree are dependent.
When producing generation $n+1$ from generation $n$, we can group the nodes into runs with a single offspring above left, and runs with a single offspring above right.
When you switch from one type of run to the other you get a node with either two or zero offspring.
We can map these runs to runs in a sequence of independent Bernoulli random variables, which have been studied in the context of various non-parametric statistics, most notably \cite{WW40}.

Huber \cite{Hub91} and Takayasu \& Takayasu \cite{TT97} have considered sums of the form $\sum_{i\in T} X_i$ where $T$ is a drainage tree arising from a diamond lattice, and the $X_i$ are i.i.d.
They use the tree to model the aggregation of charged particles.

In the case where $\be = 0$ and the drainage tree has dimension one, Equation (\ref{runoff.eqn}) is the same as the equation for the waiting time in a single server queue.
In the queuing theory literature, much use is made of the time reversed process, which satisfies the same equation as the original.
In our case reversing the process gives instead Equation (\ref{maxsum.eqn}).

Equations of type (\ref{W.eqn}) are known as Distributional Fixed Point Equations or Recursive Distributional Equations.
There is some general theory on the existence and uniqueness of solutions to such equations \cite{AB05, JOC13}, though in our case this is not an issue, as we can construct solutions using the drainage tree.

When $X$ takes values on $\{-1, 0, 1, \ldots\}$, Goldschmidt \& Pryzykucki \cite{GP17} observe that the runoff process is equivalent to a {\em parking} process, used to analyse the performance of hash tables.
They have a number of nice results for parking on a critical BGW tree with Poisson offsping numbers, on subcritical and supercritical BGW trees, and they conjecture about more general behaviour on critical trees.

\section{Proofs: the mean and right tail of $W$}\label{proofs.sec}

Since $X \in \Z$ we have $W \in \Z_+$, and we define $p_i = \Pb(W = i)$.
For convenience we also define $p_0^L = \Pb(W_LI_L = 0) = \be + \bebar p_0$ and $p_0^R = \Pb(W_RI_R = 0) = \bebar + \be p_0$.
Let $f$ be the pgf of $W$, $f(t) = \Ex t^W$, then from (\ref{W.eqn}) we have
\begin{eqnarray*}
f(t) &=& \Ex t^{I_LW_L + I_RW_R + X \vee 0} \\
&=& \Ex t^{I_LW_L + I_RW_R + 1} I_{\{X=1\}} \\
&& + \Ex I_{\{X=-1\}} I_{\{I_LW_L = 0\}} I_{\{I_RW_R = 0\}} \\
&& + \Ex t^{I_LW_L - 1} I_{\{X=-1\}} I_{\{I_LW_L > 0\}} I_{\{I_RW_R = 0\}} \\
&& + \Ex t^{I_RW_R - 1} I_{\{X=-1\}} I_{\{I_LW_L = 0\}} I_{\{I_RW_R > 0\}} \\
&& + \Ex t^{I_LW_L + I_RW_R - 1} I_{\{X=-1\}} I_{\{I_LW_L > 0\}} I_{\{I_RW_R > 0\}} \\
&=& \al (\be + (1-\be)f(t)) (1-\be + \be f(t)) t \\
&& + (1-\al) p_0^L p_0^R \\
&& + (1-\al) p_0^R (1-\be) (f(t) - p_0) t^{-1} \\
&& + (1-\al) p_0^L \be (f(t) - p_0) t^{-1} \\
&& + (1-\al) \be (1-\be) (f(t) - p_0)^2 t^{-1}  \\
&=& t^{-1} \left[ (1-\al) (1 - 2\be + 2 \be^2) p_0 (t - 1) \right. \\
&& \qquad + (1 - \al) (1 - \be) \be p_0^2 (t - 1) + (1 - \be) \be (1 + \al (t - 1)) t \\
&& + f(t) (1 - 2\be + 2\be^2) (1 + \al (t^2 - 1)) \\
&& + \left. f(t)^2 (1 - \be) \be (1 + \al (t^2 - 1)) \right]
\end{eqnarray*}

For $\be > 0$ this is just a quadratic in $f(t)$, which we can solve to give
\begin{eqnarray*}
f(t) &=& \frac {t - (\be^2+\bebar^2) (1-\al(1-t^2)) \pm \sqrt{g(t)}}
 {2\be\bebar(1-\al(1-t^2))} \\
g(t) &=& 4\albar\be^2\bebar^2 (1-\al(1-t^2)) (1-t)
\left( \left( p_0 + \frac {\be^2 + \bebar^2} {2\be\bebar} \right)^2 - h(t) \right) \\
h(t) &=& \frac {t[1 - \al(4\be\bebar + t) - \al^2(1-2\be)^2(1-t^2)]}
{4\albar\be^2\bebar^2(1 - \al(1-t^2))}
\end{eqnarray*}

To find $p_0$ and to work out which root of $g$ we use in $f$ (positive or negative), we consider how $f$ behaves at $0$ and $1$.
We will need the following result on $h$.

\begin{lem}
$h$ has a unique maximum on $[0, 1]$.

Let $t_0$ be the point at which $h$ achieves its maximum in $[0, 1]$, and define
\[
\al_c = \al_c(\be) = \frac 12 \left( 1 + \be\bebar - \sqrt{\be\bebar (2 + \be\bebar)} \right).
\]
Then for $\al > \al_c$ we have $t_0 < 1$ and $h(t_0) > h(1)$, while for $\al \leq \al_c$ we have $t_0 = 1$.
\end{lem}
{\sc Proof}
We show first that $h$ has at most one point of inflection in $[0, 1]$.
Clearly, this is equivalent to showing that $r(t) = 4\albar\be^2\bebar^2 h'(t)$ has at most one zero in $[0, 1]$.
Writing $\ga = 4\be\bebar$ we have $1-\ga = (1-2\be)^2$, so $\ga \in (0, 1]$ and
\[
r(t) = \frac {s(t)} {\left(1-\al \left(1-t^2\right)\right)^2}
= \frac {c_4 t^4 + c_2 t^2 + c_1 t + c_0} {\left(1-\al \left(1-t^2\right)\right)^2}
\]
where
\begin{eqnarray*}
c_4 &=& -\al(1-\ga) \quad \leq \quad 0 \\
c_2 &=& -\left( 2\al^2(1-\al)(1-\ga) + 2\al^2(1-\ga) + \al(1-\al) \right) \quad \leq \quad 0 \\
c_1 &=& -2\al(1-\al) \quad < \quad 0 \\
c_0 &=& (1-\al)^2 + \al(1-\ga) - \al^3(1-\ga) \quad > \quad 0.
\end{eqnarray*}
Since the denominator of $r(t)$ is strictly positive on $[0, 1]$ it is sufficient to consider the zeros of the numerator  $s(t)$.

Assuming $\be < 1/2$ so that $\ga < 1$ we can re-express the equation $s(t) = 0$ as
\[
\left( t^2 + \frac {c_2} {2c_4} \right)^2 = d_1 t + d_0
\]
for some $d_1$ and $d_0$.
Since $c_2/(2c_4) \geq 0$ this has at most two solutions.
However, $s(0) = (1-\al)(1 + \al^2(1 - \ga) - \al\ga) > 0$ and $s(t) \to -\infty$ as $t \to -\infty$, so $s$ has at least one root $< 0$, and thus at most one root in $[0, 1]$.

If $\be = 1/2$ then $\ga = 1$ and $s(t) = -\al(1-\al) (1 + t)^2 + 1 - \al$, which has roots $\pm \sqrt{1/\al} - 1$, at most one of which lie in $[0, 1]$.

Now $h(0) = 0$ and $h'(0) = c_0/(4\albar^3\be^2\bebar^2) > 0$, so $t_0 > 0$.
Since $h$ has at most one inflection point in $[0, 1]$, it follows that $t_0 = 1$ precisely when $h'(1) \geq 0$.
We have
\begin{equation}\label{hdash1.eqn}
h'(1) = \frac {1 + 4\al^2 - 4\al(1 + \be(1-\be))} {4(1-\al)\be^2(1-\be)^2}.
\end{equation}
Thus $h'(1) < 0$ iff
\[
4\al^2 - 4\al(1 + \be(1-\be)) + 1 < 0.
\]
That is, on inspecting the roots of the LHS, $h'(1) < 0$ iff
\begin{eqnarray*}
\al &>& \frac 12 \left( 1 + \be(1-\be) - \sqrt{\be(1-\be) (2 + \be(1-\be))} \right) \\
&=& \al_c(\be) \mbox{ say}.
\end{eqnarray*}
\ \hfill $\Box$
\bigskip

We can easily check that $\al_c$ is monotonic in $\be$, with $\al_c(0) = 1/2$ and $\al_c(1/2) = 1/4$.
A plot is given in Figure \ref{alpha-c.fig}.

\begin{prop}
If $\al \leq \al_c(\be)$ then $f(t)$ uses the positive root of $g(t)$ for all $t \in [0, 1]$ and
\begin{eqnarray}
p_0 &=& \frac {2\be(1-\be) - 1 + \sqrt{1 - 4\be(1-\be)\al/(1-\al)}} {2\be(1-\be)} \label{p0.eqn1} \\
\Ex W &=& \frac 1{2\be(1-\be)} \left( 1 - 2\al - \sqrt{1 - 4\al(1 - \al + \be(1-\be))} \right). \label{ExW.eqn1}
\end{eqnarray}
If $\al > \al_c(\be)$ then $t_0 = \arg\max_{t\in [0,1]} h(t) \in (0, 1)$ and $f(t)$ uses the positive root of $g(t)$ on $[0, t_0]$ and the negative root on $[t_0, 1]$, and
\begin{eqnarray}
p_0 &=& \sqrt{h(t_0)} - \frac {\be^2+\bebar^2} {2\be\bebar} \label{p0.eqn2} \\
\Ex W &=& \infty. \label{ExW.eqn2}
\end{eqnarray}
\end{prop}
{\sc Proof}
At $t = 0$ we have $h(0)=0$, $g(0) = (1-\al)^2 (\bebar^2 + \be^2 + 2\be\bebar p_0)^2$, and thus
\begin{eqnarray*}
f(0) &=& \frac {-(1-2\be+2\be^2) (1-\al) \pm \sqrt{g(0)}}
 {2(1-\be)\be(1-\al)} \\
&=& \frac {-(1-2\be+2\be^2) (1-\al) \pm (1-\al)((1-2\be+2\be^2) + 2\be(1-\be)p_0)}
 {2(1-\be)\be(1-\al)}
\end{eqnarray*}
Since $f(0) = p_0 > 0$ we must have that at $0$, $f$ uses the positive root of $g$.
Moreover, since $f$ is continuous on $[0, 1]$ and $g(0) > 0$, $f$ uses the positive root in a neighbourhood of $0$.

\begin{figure}
\begin{center}
\includegraphics[width=10cm]{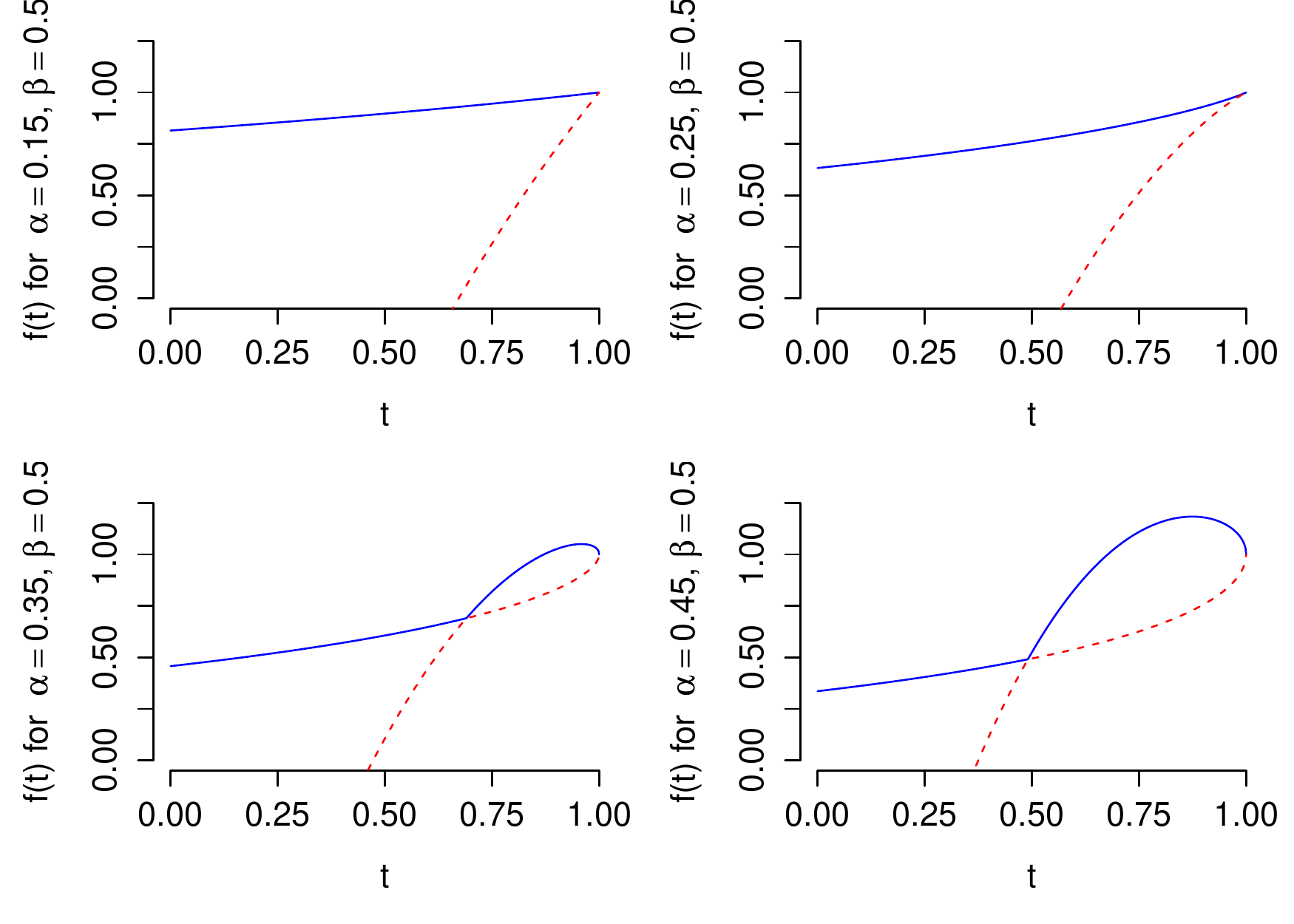}
\caption{The two branches of $f$ in the case $\be = 0.5$ and for $\al \lesseqqgtr \al_c = 0.25$.
The solid line is the branch using the positive root of $g$, and the dashed line the branch using the negative root.}\label{f-branches.fig}
\end{center}
\end{figure}

We will see that if $t_0$, the point where $h$ attains its maximum on $[0,1]$, is $< 1$, then the root used by $f$ switches at that point.
(In Figure \ref{f-branches.fig} we plot the two branches of $f$ for $\al < \al_c$, $\al = \al_c$ and $\al > \al_c$.)
First note that because $f$ is real we must have $g$ non-negative on $[0,1]$, whence
\[
h^* := \left( p_0 + \frac {\be^2 + \bebar^2} {2\be\bebar} \right)^2 \geq h(t_0) = \max_{t\in [0,1]} h(t).
\]
Now consider the behaviour of $f$ near $1$.
\begin{eqnarray*}
f(t) &=&
\frac {t - (1-2\be+2\be^2) (1-\al(1-t^2)) \pm \sqrt{g(t)}}
 {2(1-\be)\be(1-\al(1-t^2))} \\
f'(t) &=&
\frac {1 - (1-2\be+2\be^2) 2\al t \pm g'(t)/(2\sqrt{g(t)})} {2(1-\be)\be(1-\al(1-t^2))} \\
&& - \frac {t - (1-2\be+2\be^2) (1-\al(1-t^2)) \pm \sqrt{g(t)}}
 {(2(1-\be)\be(1-\al(1-t^2)))^2} 2\be(1-\be) 2\al t \\
g(t) &=&
4\be^2(1-\be)^2 (1-\al)(1-t)(1-\al(1-t^2)) (h^* - h(t)) \\
g'(t) &=&
-h'(t) 4\be^2(1-\be)^2 (1-\al)(1-t)(1-\al(1-t^2)) \\
&& + (h^* - h(t)) 4\be^2(1-\be)^2 (1-\al) ((1-t)2\al t - (1 - \al(1-t^2)))
\end{eqnarray*}
We have $f(1) = 1$, $g(1) = 0$ and
\begin{eqnarray*}
f'(1) &=& \frac {1 - (\be^2+\bebar^2)2\al \pm \lim_{t\uparrow 1} g'(t)/(2\sqrt{g(t)})} {2\be\bebar} - 2\al \\
g'(1) &=& -(h^* - h(1))4\be^2\bebar^2(1 - \al).
\end{eqnarray*}
We must have $f'(1) \geq 0$, so if $h^* > h(1)$ then $f$ must take the negative root of $g$ near 1, and $f'(1) = \infty$.

If $t_0 < 1$ then $h^* \geq h(t_0) > h(1)$, so the root of $g$ used by $f$ switches from the positive to the negative at some point in $(0,1)$.
Since $f$ is continuous on $[0, 1]$ we must have that $g(t) = 0$ at the point where the root switches.
That is, we must have $h(t) = h^*$ at this point.
It follows that the root switches at $t_0$ and that $h^* = h(t_0)$.
That is, if $\al > \al_c$ then $p_0$ solves $h^* = h(t_0)$ and $\Ex W = f'(1) = \infty$.

If $t_0 = 1$ then, as it is continuous, $f$ must use the positive root of $g$ on $[0, 1)$.
Thus we must have $h^* = h(1)$ because otherwise we would get $f'(t) < 0$ somewhere to the left of $1$.
We get $p_0$ in this case by solving $h^* = h(1)$, that is
\[
\left( p_0 + \frac {\be^2 + \bebar^2} {2\be\bebar} \right)^2
= \frac {1 - 2\al\be\bebar} {4(1-\al)\be^2\bebar^2}.
\]
To obtain $f'(1)$ in this case put $h(1) - h(t) = (1-t)(h'(1) + o(1))$, then as $t \uparrow 1$,
\begin{eqnarray*}
\frac {g'(t)}{\sqrt{g(t)}} &=&
\frac { \begin{array}{l}
2\be(1-\be)(1-\al)[-(1-\al(1-t^2))h'(t) + (1-t)(h'(1)+o(1))2\al t \\
\qquad\qquad\qquad\qquad\qquad - (1-\al(1-t^2)(h'(1) + o(1))]
\end{array} } {\sqrt{ (1-\al(1-t^2)(h'(1) + o(1)) }} \\
\lim_{t\uparrow 1} \frac {g'(t)}{\sqrt{g(t)}} &=&
-4\be(1-\be)\sqrt{(1-\al)h'(1)}
\end{eqnarray*}
Thus, plugging in $h'(1)$ (see Equation (\ref{hdash1.eqn})), we have
\begin{equation}\label{ExW:eqn}
\Ex W
= f'(1)
= \frac 1{2\be(1-\be)} \left( 1 - 2\al - \sqrt{1 - 4\al(1 - \al + \be(1-\be))} \right).
\end{equation}
Plots of $\Ex W$ for various $\al$ and $\be$ are given in Figure \ref{ExW.fig}.
\ \hfill $\Box$
\bigskip

\begin{rem}
Note that the expression for $\Ex W$ simplifies for $\al = \al_c$, giving
\begin{equation}\label{ExWcrit.eqn}
\frac 12 \left( \sqrt{ 1 + \frac 2{\be(1-\be)} } - 1 \right).
\end{equation}
Also note that for $\be=1/2$ the expression for $p_0$ has a simple form for all $\al$ (here $\al_c(1/2) = 1/4$):
\[
p_0 = \left\{\begin{array}{ll}
2 \sqrt{ \frac {1-2\al}{1-\al} } - 1 & \mbox{ for } 0 \leq \al \leq 1/4 \\
\sqrt{ \frac 2{\sqrt{\al} + \al} } - 1 & \mbox{ for } 1/4 \leq \al \leq 1.
\end{array} \right.
\]
For $\be = 1/2$ and $\al \leq 1/4$ we also have $\Ex W = 2(1 - 2\al - \sqrt{(1 - \al)(1 - 4\al)})$.
\end{rem}

\begin{prop}\label{main.prop}
Let $F$ be the cdf of $W$ then, as $x \to \infty$,
\begin{equation}
1 - F(x) \sim
\left\{\begin{array}{ll}
\sqrt{ \frac {-h''(1)(1-\al)}{8\pi}} x^{-3/2} & \al = \al_c \\
\sqrt{ \frac {(h(t_0) - h(1))(1-\al)}\pi } x^{-1/2} & \al > \al_c
\end{array}\right.
\end{equation}
For $\al < \al_c$ $W$ has all positive moments finite.
\end{prop}
{\sc Proof}
We work with Laplace-Stieltjes transforms.
Let $\Fhat$ be the L-S transform of $F$, so $\Fhat(s) = f(e^{-s})$.
We are interested in the behaviour of $\Fhat(s)$ near $0$, that is, the behaviour of $f$ near $1$.
We will write $P_k$ to indicate a generic polynomial whose smallest non-zero term is order $k$,
possibly of infinite order, but convergent in a neighbourhood of $0$.

The behaviour of $f$ at 1 depends on the behaviour of $g$, which in turn depends on the term $h(t_0) - h(t)$.
If $\al < \al_c$ then $t_0 = 1$ and $h'(1) > 0$, so $h(t_0) - h(t) = h'(1)(1 - t) + P_2(1-t)$.
If $\al = \al_c$ then $t_0 = 1$ and $h'(1) = 0$, so $h(t_0) - h(t) = -\half h''(1) (1 - t)^2 + P_3(1-t)$.
If $\al > \al_c$ then $t_0 < 1$ and $h(1) < h(t_0)$, so $h(t_0) - h(t) = h(t_0) - h(1) + P_1(1- t)$.
We take each case in turn.
In each case we use the fact that near 0, $1 - e^{-s} = s + P_2(s)$.

For $\al < \al_c$, $f$ uses the positive root of $g$ near 1, and
\[
h(1) - h(e^{-s}) = (1-e^{-s}) h'(1) + P_2(1-e^{-s}) = s h'(1) + P_2(s).
\]
This gives
\[
\Fhat(s) = P_0(s) + s\sqrt{P_0(s)}
\]
which has a convergent Taylor series expansion in a neighbourhood of $s=0$.
Thus $\Fhat$ has all its derivatives finite at $0$, so $W$ has all positive moments finite.

In the case $\al = \al_c$ $f$ again uses the positive root of $g$ near 1, and
\[
h(1) - h(e^{-s})
= -\half(1 - e^{-s})^2h''(1) + P_3(1-e^{-s})
= -\half s^2 h''(1) + P_3(s),
\]
Thus
\[
\Fhat(s) = \frac {2\be\bebar + (2\al(1 - 2\be\bebar) - 1)s + O(s^2) + s^{3/2} \be\bebar \sqrt{-2(1-\al) h''(1) + O(s)}} {2\be\bebar(1 - 2\al s) + O(s^2)}.
\]
Writing $\mu$ for $\Ex W$ and plugging in our expression for $\al = \al_c$ we get, for $s \downarrow 0$,
\begin{eqnarray*}
\Fhat(s) - 1 + \mu s
&=& s^{3/2} \sqrt{-\half (1-\al) h''(1) + O(s)} + O(s^2) \\
&=& s^{3/2} l(1/s)
\end{eqnarray*}
where $l$ is slowly varying at infinity.
Standard Tauberian theory now tells us (see Bingham, Goldie \& Teugels \cite{BGT87} Theorem 8.1.6) that as $x \to \infty$
\begin{eqnarray*}
1 - F(x) &\sim& l(x) \frac {\Ga(3/2)}{\Ga(1/2)^2} x^{-3/2} \\
&\sim& \sqrt{ \frac {-h''(1)(1-\al)}{8\pi}} x^{-3/2}.
\end{eqnarray*}

In the case $\al > \al_c$ $f$ uses the negative root of $g$ near 1, and we get
\begin{eqnarray*}
f(t) &=& 1 - \sqrt{1-t} \sqrt{(h(t_0) - h(1))(1-\al) + O(1-t)} + O(1-t)\\
f(e^{-s}) - 1 &=& - s^{1/2} \sqrt{(h(t_0) - h(1))(1-\al) + O(s)} + O(s).
\end{eqnarray*}
That is, $\Fhat(s) - 1 = s^{1/2} l(1/s)$ where $l$ is slowly varying, so applying our Tauberian theorem we see that as $x\to\infty$
\[
1 - F(x) \sim \sqrt{ \frac {(h(t_0) - h(1))(1-\al)}\pi } x^{-1/2}.
\]
\ \hfill $\Box$
\bigskip

\subsection{Case: $\be=0$}\label{be0.sec}

When $\be = 0$ and $\al \in[0, 1/2)$ we get $f(t) = \albar p_0/(\albar - \al t)$, from which it follows that $p_0 = (1-2\al)/(1-\al)$ (since $f(1) = 1$), and thus that $W \sim \mbox{geom}((1-2\al)/(1-\al))$ and $\Ex W = \al/(1-2\al)$.
For $\al \in [1/2, 1]$ we have that $p_0 = 0$ and $W = \infty$ almost surely.

\section{How much of the tree contributes to root runoff?}

In this section we look at how much of the tree is contributing to the runoff at the root; our results are summarised in Proposition \ref{contribution.prop}.
Recall that $W = (W_LI_L + W_RI_R + X) \vee 0 = W_LI_L + W_RI_R + Y$, where $Y$ is the nett contribution to runoff from our cell.
If $X = 1$ then $Y = 1$.
If $W_L I_L + W_R I_R = 0$ and $X = -1$ then $Y = 0$.
If $W_L I_L + W_R I_R > 0$ and $X = -1$ then $Y = -1$.
Thus $Y$ has distribution
\begin{eqnarray*}
\Pb(Y = y) &=&
\left\{\begin{array}{ll}
\al, & \quad y = 1 \\
(1-\al) [\be + (1-\be)p_0] [1-\be + \be p_0], & \quad y = 0 \\
(1-\al) (1 - [\be + (1-\be)p_0] [1-\be + \be p_0]), & \quad y = -1
\end{array}\right. \\
\Ex Y &=& 2\al - 1 + (1-\al)(\be + (1-\be)p_0)(1-\be + \be p_0).
\end{eqnarray*}

It is easy to check that
\[
\Ex Y \left\{\begin{array}{ll}
= 0 & \al \leq \al_c \\
> 0 & \al > \al_c
\end{array}\right.
\]
In particular for $\be = 1/2$ and $\al > \al_c = 1/4$ we get $\Ex Y = {(1 + \sqrt{\al})(2\sqrt{\al} - 1)^2}/{(2\sqrt{\al})} > 0$.
Thus when $\al > \al_c$ the nett contribution at each node has positive mean, and we expect most of the tree to be contributing to runoff at the root.

Note that if we {\em know} $\Ex W < \infty$ then we get $\Ex Y = 0$ from (\ref{Y.eqn}), which then gives us the same equation for $p_0$ as in (\ref{p0.eqn1}) for $\al \leq \al_c$.

\subsection{Size of the tree}\label{size_tree.sec}

Further evidence that most of the tree contributes to the root-runoff when $\al > \al_c$ comes from comparing the right tails of $W$ and $N_T$, the total number of nodes in the tree.
It is known that for a (sub)critical GW process with offspring distribution $\xi$, the total progeny $N_T$ has the same law as $T_1$, the first time to hit $-1$ for a r.w.\ with steps distributed as $\xi - 1$, started at 0 \cite{LeG05}.
Let $\chi_i$ be i.i.d.\ distributed as $\xi - 1$ and put $S_n = \sum_{i=1}^n \chi_i$.
Since our r.w. is left-continuous we have \cite{HK08}
\[
\Pb(T_1 = n) = \frac 1n \Pb(S_n = -1)
\]

In our case we have $S_n \sim M_1 - M_3$ where $(M_1, M_2, M_3) \sim \mbox{multinomial}(n, (\be\bebar, \be^2+\bebar^2, \be\bebar))$.
Moreover, for large $n$, writing $\rho = \be\bebar$,
\[
\left( \begin{array}{c} M_1 \\ M_2 \\ M_3 \end{array} \right)
\approx N\left( 
\left( \begin{array}{c} n\rho \\ n(1-2\rho) \\ n\rho \end{array} \right),
\left( \begin{array}{ccc}  n\rho(1-\rho) & -n\rho(1-2\rho) & -n\rho^2 \\
-n\rho(1-2\rho) & n2\rho(1-2\rho) & -n\rho(1-2\rho) \\
-n\rho^2 & -n\rho(1-2\rho) & n\rho(1-\rho) \end{array} \right)
\right)
\]
Thus $M_1 - M_3 \approx N(0,2n\be\bebar)$ and (using the usual continuity correction)
\begin{eqnarray*}
\Pb(M_1 - M_3 = -1) &\approx& \Phi\left( \frac 3{2\sqrt{2n\be\bebar}} \right) - \Phi\left( \frac 1{2\sqrt{2n\be\bebar}} \right) \\
&\approx& \phi(0) \frac 1{\sqrt{2n\be\bebar}}
\ =\ \frac 1{2\sqrt{\pi n\be\bebar}}.
\end{eqnarray*}
That is
\[
\Pb(N_T = n) \approx n^{-3/2} \frac 1{2\sqrt{\pi\be\bebar}}.
\]
Thus, if some fixed percentage of the whole tree was contributing to the root-runoff, we would expect $1 - F(x) = O(x^{-1/2})$, which is indeed the case (from (\ref{main.eqn})).

\subsection{Runoff down the spine}

We can show how the drainage tree contributes to runoff by considering its {\em spine} \cite{LPP95, GK99}.

Condition the tree to be of height $n$, then consider the left-most line of descent of length $n$.
At any fixed depth along this line of descent (the spine), the offspring distribution converges to the size-biased offspring distribution as $n \to \infty$.
Subtrees attached to the spine grow like the original tree but with limited height: at generation $k$ subtrees growing to the left are conditioned to have height at most $n - k - 1$, while subtrees growing to the right are conditioned to have height at most $n - k$.
So, fixing $k$ and sending $n \to \infty$, the runoff coming from a subtree attached to the spine at generation $k$ will have a distribution tending to $W$.
In our case the size biased distribution is 1 with probability $1 - 2\be\bebar$ and 2 with probability $2\be\bebar$.
That is, spinal nodes have can have at most one subtree attached, with probability tending to $2\be\bebar$ as $n\to\infty$.

We now consider the runoff process on the spine.
At each point on the spine we have a point contribution, distributed as $X$, and with some positive probability, runoff generated by a subtree.
If $\al > \al_c$ then, fixing the generation $k$ and sending the tree height $n \to \infty$, the mean runoff from a subtree will tend to infinity.
It follows that for any $\de > 0$ we can choose an $m$ such that for all $n \geq m$, the runoff process down the spine is bounded below by a random walk with positive drift, at least for the bottom $n(1-\de)$ nodes.
Thus there will be a positive probability that the runoff will be strictly positive all the way down the bottom $100(1-\de)\%$ of the spine.

If $\al < \al_c$ then the mean runoff from subtrees must be finite, so as $n \to \infty$ the mean runoff generated at each point on the spine will tend to
\[
\de := \Ex X + 2\be\bebar \Ex W = 2\al - 1 + 2\be\bebar \Ex W.
\]
But from (\ref{ExWcrit.eqn}), for $\al < \al_c = (1 + \be\bebar - \sqrt{\be\bebar(2 + \be\bebar)})/2$ we have that $\Ex W < (\sqrt{(2+\be\bebar)/(\be\bebar)} - 1)/2$, so
\[
\de < 1 + \be\bebar - \sqrt{\be\bebar(2 + \be\bebar)} - 1 + \sqrt{(2+\be\bebar)\be\bebar} - \be\bebar = 0.
\]
Thus if $\al < \al_c$ then runoff down the spine behaves like the waiting time process in a stable single server queue.
In this case the expected queue size is finite, which translates as saying the number of nodes on the spine that contribute to the runoff at the root, has finite mean.

\section{Left continuous $X$}\label{left_cts.sec}

The approach taken in Section \ref{proofs.sec} for $X \in \{-1, 1\}$ can be largely extended to left continuous variables.
Suppose that $X \in \{-1, 0, 1, \ldots \}$.
Let $\al = \Pb( X \geq 0)$, $m = \Ex X$ and put $\eta(t) = \Ex t^{X+1}$ (a proper pgf).
Let $f(t) = \Ex t^W$ as before, then conditioning on $\{ I_LW_L = 0\}$, $\{ I_RW_R = 0\}$ and $\{ X = -1\}$, we get
\begin{eqnarray*}
f(t)
&=& f(t)^2 \be\bebar \Ex t^X + f(t) (\be^2+\bebar^2) \Ex t^X \\
&& + \be\bebar \Ex t^X + (1-\al) (\be+\bebar p_0) (\bebar + \be p_0) (1 - t^{-1})
\end{eqnarray*}
Solving for $f$ we get
\begin{eqnarray*}
f(t) &=& \frac {t - (\be^2+\bebar^2) \eta(t) \pm \sqrt{g(t)}} {2\be\bebar\eta(t)} \\
g(t)
&=& 4\be^2\bebar^2(1-\al)(1-t)\eta(t) \left[ \left( p_0 + \frac {\be^2+\bebar^2}{2\be\bebar} \right)^2 - h(t) \right] \\
h(t) &=& \frac 1{4\be^2\bebar^2(1-\al)\eta(t)} \left( (1-2\be)^2 \eta(t) \frac {1-\eta(t)}{1-t} + (1-\al)(1-2\be)^2 \eta(t) \right. \\
&& \hspace{4cm} \left. -2(\be^2+\bebar^2) \eta(t) - \frac {1-\eta(t)}{1-t} + 1 + t \right)
\end{eqnarray*}

To keep the algebra manageable, we will restrict ourselves to the case $\be = 1/2$ for the rest of this section.
We will also assume that $X$ has finite mean and variance.
When $\be=1/2$ we have
\begin{eqnarray*}
f(t) &=& \frac {2t - \eta(t) \pm 2\sqrt{g(t)}} {\eta(t)} \\
g(t) &=& \frac 14 (1-\al) (1-t) \eta(t) \left( (p_0+1)^2 - h(t) \right) \\
h(t)
&=& \frac 4{(1-\al)} \frac {t(\eta(t) - t)}{(1-t)\eta(t)}
\end{eqnarray*}

\begin{prop}\label{left_cts_X.prop1}
Assume that $X \in \{-1, 0, 1, \ldots \}$ has finite mean and variance, and that $h$ has a unique max in $[0, 1]$, occurring at $t_0 > 0$.
If $t_0 = 1$ then $f$ takes the positive root of $g$ on $[0, 1]$ and, writing $m = \Ex X$ (necessarily negative in this case),
\begin{eqnarray*}
p_0 &=& 2\sqrt{\frac {-m}{1-\al}} - 1 \\
\Ex W &=& -2m - \sqrt{2(-m(1-m) - \Var X)}.
\end{eqnarray*}
If $t_0 \in (0, 1)$ then $f$ takes the positive root of $g$ on $[0, t_0]$ and the negative root on $[t_0, 1]$, and
\begin{eqnarray*}
p_0 &=& \sqrt{h(t_0)} - 1 \\
\Ex W &=& \infty.
\end{eqnarray*}
\end{prop}
{\bf Proof}
Consider first the behaviour of $f$ at $0$.
We have
\begin{eqnarray*}
\eta(0) &=& 1 - \al \\
h(0) &=& 0 \\
g(0) &=& \frac 14 (1-\al)^2 \left( p_0 + 1 \right)^2 \\
f(0) &=& -1 \pm (p_0 + 1)
\end{eqnarray*}
Thus, since $f(0) = p_0$, near $0$ $f$ must take the +ve root of $g$.

Now consider the case $t_0 = 1$.
Since $g(t) \geq 0$ we have
\[
(p_0 + 1)^2 \geq \max_{0\leq t \leq 1}  h(t) = h(1).
\]
Moreover, assuming $m = \Ex X < \infty$, we have that near $t = 1$, $\eta(t) = 1 - (1-t)(m+1) + o(1-t)$ and so
\begin{eqnarray*}
h(1) &=& \frac {-4m}{1-\al} \\
p_0 &\geq& \sqrt{ \frac {-4m}{1-\al} } - 1
\end{eqnarray*}

By assumption $h$ has a unique maximum at $t_0 = 1$, so $g(t) > 0$ for all $t \in [0, 1)$.
Thus since $f$ is continuous it uses the positive root of $g$ for all $t \in [0, 1]$, whence near $t = 1$ we have
\[
f'(t) = \frac {2 - \eta'(t) + g'(t)/\sqrt{g(t)}} {\eta(t)}
- \frac {2t - \eta(t) + 2\sqrt{g(t)}}{\eta(t)^2} \eta'(t).
\]
As $t \uparrow 1$ we have $\eta(t) \uparrow 1$, $\eta'(t) \uparrow m + 1$, $g(t) \downarrow 0$, so
\[
\lim_{t \uparrow 1} f'(t) = -2m + \lim_{t \uparrow 1} \frac {g'(t)}{\sqrt{g(t)}}.
\]
Now, provided $h'(1)$ is finite, near $t = 1$ we have
\begin{eqnarray*}
g'(t) &=& \frac {1-\al}4 \left[ (-\eta(t) + (1-t)\eta'(t)) ((p_0+1)^2 - h(t)) - (1-t) \eta(t) h'(t) \right] \\
&=& -\frac {(1-\al) \eta(t)}4 ((p_0+1)^2 - h(t)) + O((1-t))
\end{eqnarray*}
Thus, since $\lim_{t \uparrow 1} f'(t) \geq 0$, we must have
\[
(p_0+1)^2 \leq h(1).
\]
That is, if $h$ achieves its maximum on $[0, 1]$ at 1 (and nowhere else), and $h'(1)$ is finite, then $(p_0+1)^2 = h(1)$.
That is, $p_0 = \sqrt{ -4m/(1-\al) } - 1$.

Assuming $\Var X < \infty$, near $t=1$ we have $\eta(t) = 1 - (1-t)(m+1) + \half (1-t)^2 \eta''(1) + o((1-t)^2)$.
Thus
\begin{eqnarray*}
h(1) - h(t)
&=& \frac 4{1-\al} (1 - t) \left( m^2 - \half \eta''(1) \right)  + o(1-t) \\
&=& (1-t) \frac 2{1-\al} (-m(1-m) - \Var X) + o(1-t)
\end{eqnarray*}
So $h'(1) = 2(-m(1-m) - \Var X)/(1-\al)$ is finite as required.

Plugging our value for $p_0$ into the expression for $f'(1) = \Ex W$ we get
\begin{eqnarray*}
g(t) &=& \frac {1-\al}4 h'(1) (1-t)^2 + o((1-t)^2) \\
g'(t) &=& - \frac {1-\al}{2} h'(1) (1-t) + o(1-t) \\
f'(1) &=& -2m + \lim_{t\uparrow 1} \frac {g'(t)}{\sqrt{g(t)}} \\
&=& -2m - \sqrt{(1-\al)h'(1)} \\
&=& -2m - \sqrt{2(m(m-1) - \Var X)}
\end{eqnarray*}

Now consider the case $t_0 \in (0, 1)$.
We have $h(t_0) > h(1)$ so near $t = 1$
\[
g'(t) = - \frac {1-\al}4 ((p_0+1)^2 - h(1)) + O((1-t)) \ < \ 0
\]
Thus, since $\lim_{t\uparrow 1} f'(t) \geq 0$, we must have that near $t=1$ $f$ takes the negative root of $g$, whence
\[
\Ex W
\ =\ f'(1)
\ =\ -2m - \lim_{t\uparrow 1} \frac {g'(t)}{\sqrt{g(t)}}
\ =\ \infty
\]
Also, because $f$ is continuous, we must have $g(t)=0$ at the point where the root switches.
That is,
\[
p_0 = \sqrt{h(t_0)} - 1.
\]
\ \hfill $\Box$
\bigskip

\begin{prop}
Suppose that $\eta$ exists in a neighbourhood of 1, and that $h$ has a unique max in $[0, 1]$, occurring at $t_0 > 0$.
If $t_0 = 1$ and $\Var X < -m(1-m)$ then $W$ has all positive moments finite (subcritical case). 

If $t_0 = 1$ and $\Var X = -m(1-m)$ (critical case) then, as $x \to \infty$, 
\[
1 - F(x) \sim \sqrt{ - \frac {h''(1)(1-\al)}{8\pi}} x^{-3/2}
\]

If $t_0 < 1$ (supercritical case) then, as $x \to \infty$,
\[
1 - F(x) \sim x^{-1/2} \sqrt{ \frac {(h_{\max} - h(1))(1-\al)}\pi }
\]
\end{prop}
{\sc Proof}
From our assumption on $\eta$, $X$ has all positive moments finite.
We recall from the proof of Proposition \ref{left_cts_X.prop1} that $h'(1) = 2(-m(1-m) - \Var X)/(1-\al)$.
As before, let $\hat{F}$ be the L-S transform of $F$.
Our proof follows that of Proposition \ref{main.prop}.

Consider first the case $t_0 = 1$ and $h'(1) > 0$, that is $\Var X < -m(1-m)$.
We have
\begin{eqnarray*}
g(t) &=& \frac {1-\al}4 h'(1) (1-t)^2 + o((1-t)^2) \\
f(t) &=& \frac {2t - \eta(t) + 2\sqrt{g(t)}}{\eta(t)} \\
&=& (1 + (1-t)(m - 1) + (1-t)\sqrt{(1-\al)h'(1) + P_1(1-t)}) \\
&& \qquad \times (1 + (1-t)(m+1) + P_2(1-t)) \\
\Fhat(s) &=& (1 + s(m-1) + s\sqrt{(1-\al)h'(1) + P_1(s)} + P_2(s)) \\
&& \qquad \times (1 + s(m+1) + P_2(s)) \\
&=& P_0(s) + s\sqrt{P_0(s)}
\end{eqnarray*}
This has a convergent Taylor series expansion about 0, so $W$ has all positive moments finite.

Next we take the case $t_0 = 1$ and $h'(1) = 0$, that is $\Var X = -m(1-m)$.
We have, from Tauberian theory,
\begin{eqnarray*}
g(t) &=& -\frac {1-\al}8 h''(1) (1-t)^3 + o((1-t)^3) \\
g'(t) &=& + \frac {1-\al}4 h''(1) (1-t)^2 + o((1-t)^2) \\
\Ex W \ =\ f'(1) &=& -2m + \lim_{t\uparrow 1} \frac {g'(t)}{\sqrt{g(t)}} \\
&=& -2m \\
f(t) &=& \frac {2t - \eta(t) + 2\sqrt{g(t)}}{\eta(t)} \\
&=& (1 + (1-t)(m - 1) + (1-t)^{3/2}\sqrt{-(1-\al)h''(1)/2 + P_1(1-t)}) \\
&& \qquad \times (1 + (1-t)(m+1) + P_2(1-t)) \\
\Fhat(s) &=& (1 + s(m-1) + s^{3/2}\sqrt{-(1-\al)h''(1)/2 + P_1(s)} + P_2(s)) \\
&& \qquad \times (1 + s(m+1) + P_2(s)) \\
\Fhat(s) - 1 -2ms &=& s^{3/2}\sqrt{-(1-\al)h''(1)/2 + P_1(s)} + P_2(s) \\
1 - F(x) &\sim& \sqrt{ - \frac {h''(1)(1-\al)}{8\pi}} x^{-3/2}
\end{eqnarray*}

Lastly we take the case case $t_0 < 1$.
Again using Tauberian theory we have
\begin{eqnarray*}
g(t) &=& \frac 14 (1-\al) (1-t) (h_{\max} - h(1)) + o(1-t) \\
f(t) &=& \frac {2t - \eta(t) - 2\sqrt{g(t)}}{\eta(t)} \\
&=& (1 + (1-t)(m - 1) - (1-t)^{1/2}\sqrt{(1-\al)(h_{\max} - h(1)) + P_1(1-t)}) \\
&& \qquad \times (1 + (1-t)(m+1) + P_2(1-t)) \\
\Fhat(s) &=& (1 + s(m-1) - s^{1/2}\sqrt{(1-\al)(h_{\max} - h(1)) + P_1(s)} + P_2(s)) \\
&& \qquad \times (1 + s(m+1) + P_2(s)) \\
\Fhat(s) - 1 &=& - s^{1/2}\sqrt{(1-\al)(h_{\max} - h(1)) + P_1(s)} + P_1(s) \\
1 - F(x) &\sim& x^{-1/2} \sqrt{ \frac {(h_{\max} - h(1))(1-\al)}\pi }
\end{eqnarray*}
\ \hfill $\Box$
\bigskip

\begin{rem}
Note that if $h$ has at most one point of inflection in $[0, 1]$, then we can determine supercritical/critical/subcritical behaviour from $h'(1)$.
Specifically, if $h'(1) < 0$ then $t_0 < 1$ and the runoff is supercritical, if $h'(1) = 0$ then $t_0 = 1$ and the runoff is critical, and if $h'(1) > 0$ then $t_0 = 1$ and the runoff is subcritical.
\end{rem}

\begin{figure}
\begin{center}
\includegraphics[width=8cm]{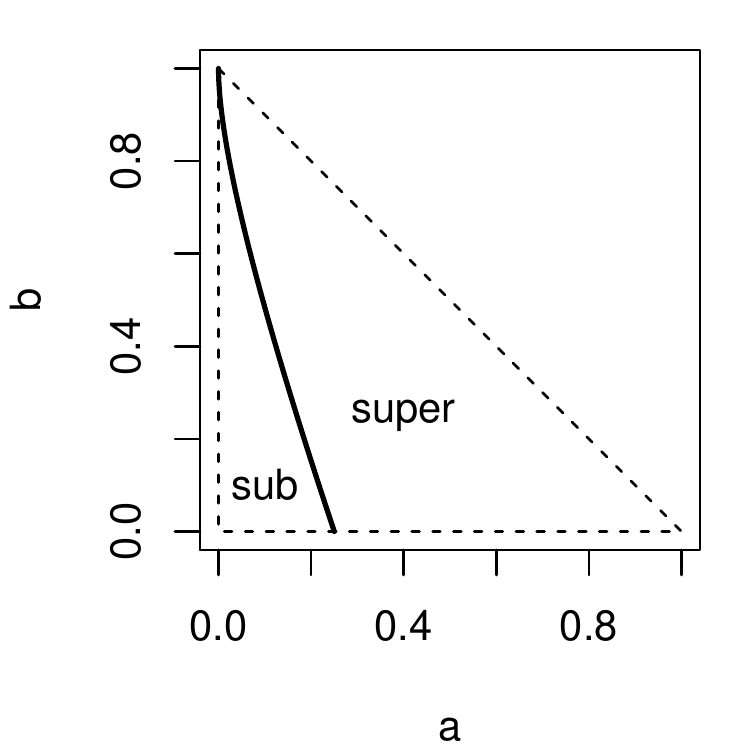}
\caption{Subcritical/critical/supercritical reg\'{e}mes for Example \ref{ex1}.
The critical region is given by the solid curve.} \label{example.fig}
\end{center}
\end{figure}

\begin{exa}\label{ex1}
Suppose that
\[
X = \left\{\begin{array}{rl}
1 & \mbox{ w.p. } a \\ 0 & \mbox{ w.p. } b \\ -1 & \mbox{ w.p. } c = 1 - a - b
\end{array}\right.
\]
then we have
\begin{eqnarray*}
h(t) &=& \frac {4 t (1 - a - b - at)} {(1 - a - b) (1 - b(1-t) - a(1-t^2))} \\
h'(1) &=& \frac {4 ((1-b)^2 + 4a^2 - a(5 - 4b))} {1 - a - b}
\end{eqnarray*}

In this case $h'$ has at most one root in $[0, 1]$, so the runoff is subcritical/critical/supercritical according to $h'(1) >/=/< 0$.
In Figure \ref{example.fig} we plot these regions as functions of $a$ and $b$.
\end{exa}

\subsection*{Acknowledgements}
Thanks to Ben Hambly for the helpful discussions and to Christina Goldschmidt for sharing her related work.


\end{document}